\documentclass{amsart}

\usepackage{amsmath}
\usepackage{amsthm}
\usepackage{amssymb}

\theoremstyle{plain}   
\begingroup 

\newtheorem{theorem}{Theorem}[section]   
\newtheorem{corollary}[theorem]{Corollary}     
\newtheorem{lemma}[theorem]{Lemma}         
\endgroup

\theoremstyle{definition}

\theoremstyle{remark}
\newtheorem{remark}[theorem]{Remark}        

\numberwithin{equation}{section}

\newcommand{\R}{{\mathbb R}}
\newcommand{\N}{{\mathbb N}}

\newcommand{\e}{\varepsilon}

\newcommand{\diam}{\operatorname{diam}}

\newcommand{\Hau}{{\mathcal{H}}}

\newcommand{\mcI}{{\mathcal{I}}}

\begin{document}	

\title{Absolutely continuous functions with values in metric spaces}

\author{Jakub Duda} 

\thanks{The author was supported in part by ISF}

\address{
Department of Mathematics, Weizmann Institute
of Science, Rehovot 76100, Israel}
\email{jakub.duda@weizmann.ac.il}

\date{December 14, 2005}




\begin{abstract} 
We present a general theory of absolutely continuous
paths with values in metric spaces using the notion of
metric derivatives. Among other results, we prove analogues of
the Banach-Zarecki and Vall\'ee Poussin theorems.
\end{abstract}

\maketitle

\section{Introduction}

In a nice expository article, Varberg~\cite{Var} 
outlined an elegant approach towards the theory of real-valued absolutely continuous functions.
In the present note, we will be interested in maps $f:[a,b]\to (M,\rho)$,
where $(M,\rho)$ is a metric space. We will see that a significant
part of the theory carries over to this (very general) situation. As
we have the following:
\begin{equation}\label{embed}
\begin{split}
&\text{ every metric space } (M,\rho)
\text{ can be embedded }\\
&\text{ into a suitable Banach space }\ell_\infty(\Gamma)
\text{ for some }\Gamma
\end{split}
\end{equation}
(see e.g.~\cite[Lemma~1.1]{BL}), we could without
any loss of generality work with Banach spaces only.
\par
The main obstacle in dealing with metric spaces (or arbitrary Banach spaces)
is the absence of the Radon-Nikod\'ym property and the resulting non-existence
of derivatives. Thus, instead of the ``usual'' derivative, we have to employ
the notion of a ``metric derivative'' (which was introduced by Kirchheim in~\cite{K}).
We will need some results about this notion from~\cite{D}.
\par
Let $(M,\rho)$ be a metric space, and let $f:[a,b]\to(M,\rho)$. We say that
{\em $f$~is absolutely continuous}, provided for each $\e>0$ there exists
a $\delta>0$ such that whenever $[a_1,b_1],\dots, [a_k,b_k]$ is
a sequence of non-overlapping intervals in $[a,b]$ with
$\sum^k_{i=1} (b_i-a_i)<\delta$, then
\[\sum^k_{i=1} \rho\big(f(b_i),f(a_i)\big)<\e.\]
It easily follows that absolutely continuous functions are continuous.
\par
This paper is organized in the following way. In the second section,
we present the basic definitions and establish some auxiliary results.
In the third section,
we present the theory of absolutely continuous functions with
values in metric spaces. For example, we prove a version of the
Banach-Zarecki theorem in this context -- see Theorem~\ref{banzar}
(which was recently proved by L.~Zaj\'\i\v{c}ek
and the author in~\cite{DZ}). The current proof is different from
the one in~\cite{DZ} -- it does not use the theorem of Luzin,
but rather a generalization of ideas due to Varberg~\cite{Var}.
Among other results, we also show a version of Vall\'ee Poussin's
theorem (see Theorem~\ref{pothm})
which characterizes the situation when a composition
of two absolutely continuous functions is again absolutely continuous.

\section{Preliminary results}

By $m$ we will denote the Lebesgue measure on $\R$.
For each function $f:[a,b]\to(M,\rho)$, and for $x\in[a,b]$ we can define the variation
\[v_f(x)=\bigvee^x_a f=\sup_D \sum^{n(D)-1}_{i=0} \rho(f(x_{i}),f(x_{i+1})),\]
where the supremum is taken over all partitions $D$ of $[a,x]$ 
($D$ is a {\em partition of $[a,x]$} provided $D=\{a=x_0<x_1<\dots<x_n=b\}$,
and $n=n(D)=\# D-1$). We say that {\em $f$ has bounded variation},
provided $\bigvee^b_a f<\infty$. It is easy to see that every absolutely
continuous function has bounded variation.
\par
We will need the notion of the~``metric derivative''. Let $f:[a,b]\to(M,\rho)$. 
For $x\in[a,b]$, we define 
\[ md(f,x):=\lim_{\substack{t\to 0\\x+t\in[a,b]}} \frac{\rho\big(f(x+t),f(x)\big)}{|t|}.\]
Following~\cite{K}, we say that {\em $f$ is metrically differentiable at $x$}, provided
$md(f,x)$ exists and 
\[ \rho(f(y),f(z))-md(f,x)|y-z|=o(|x-y|+|x-z|),\ \text{ when }(y,z)\to(x,x).\]
The following is an easy consequence of~\cite[Theorem~2.6]{D}:

\begin{theorem}\label{Dthm} 
Let $f:[a,b]\to(M,\rho)$ be arbitrary. Then the 
following hold.
\begin{enumerate}
\item If $S(f):=\{x\in[a,b]:\limsup_{t\to0} |t|^{-1}\rho(f(x+t),f(x))<\infty\}$,
then there is $N\subset[a,b]$ with $m(N)=0$ such that
$f$ is metrically differentiable at all $x\in S(f)\setminus N$.
\item If $f$ has bounded variation,
then $f$ is metrically differentiable at almost all $x\in[a,b]$.

\end{enumerate}
\end{theorem}

\begin{proof}
Part~(i) is just a restatement of~\cite[Theorem~2.6]{D}.
To prove part~(ii), note that $v_f$ is differentiable 
almost everywhere in~$[a,b]$. We easily see that at each such point
we have $\limsup_{t\to0} |t|^{-1}\rho(f(x+t),f(x))<\infty$. Thus part~(i) implies that $f$ is 
metrically differentiable
at almost each $x\in[a,b]$.
\end{proof}

We will need the following simple lemma.

\begin{lemma}\label{diflem}
Let $(M,\rho)$ be a metric space, $f:[c,d]\to M$, $g:[a,b]\to[c,d]$,
$x\in[a,b]$ be such that $g'(x)\neq0$ and $md(f\circ g,x)$ exists.
Then $md(f,g(x))$ exists.
\end{lemma}

\begin{proof}
Denote $\eta=g'(x)$. By the differentiability of $g$ at $x$, we 
have
\[ g(x+h)-g(x)-\eta\,h=o(h),\ \text{ when }h\to0.\]
Thus, we can choose $\delta>0$ such that $g(y)\neq g(x)$ for $|x-y|<\delta$,
and for each $|h|<\delta$ there exists $h'\in\R$ such that $g(x+h')=g(x)+\eta\,h$.
It is easy to see that $h\to0$ if and only if $h'\to0$. 
We have 
\[g(x+h')=g(x)+\eta\,h'+o(h')=g(x)+\eta\,h,\]
and thus $h/h'\to1$ when $h\to0$.
Now, 
\[ \frac{\rho\big(f(g(x)+\eta\,h),f(g(x))\big)}{\eta\,h}
=\frac{\rho\big(f(g(x+h')),f(g(x))\big)}{\eta\,h'}\cdot\frac{h'}{h}\to\frac{md(f\circ g,x)}{\eta},\]
when $h\to0$.
Thus $md(f,g(x))$ exists.
\end{proof}

Let $(M,\rho)$ be a metric space, and $A\subset M$. 
We define the {\em Hausdorff measure $\Hau^1(A)$} as 
$\lim_{\delta\to0} \Hau^1_\delta(A)$,
where
\[ \Hau^1_\delta(A):=\inf\big\{\sum\diam(A_i):A\subset\bigcup_i A_i\text{ with }\diam(A_i)<\delta\ \forall i\big\}.\]
for $\delta>0$. It is well known (see e.g.\ \cite{F}) that $\Hau^1$ is a Borel measure on~$M$.
\par
The following is a ``metric'' version of Varberg's ``Fundamental Lemma''
(see~\cite[p.~832]{Var}).

\begin{lemma}\label{funlem}
Let $f:[a,b]\to (M,\rho)$ be a function, let $E$ be the
set of all $x\in[a,b]$ where $md(f,x)$ exists and satisfies $md(f,x)\leq K$.
Then
\begin{equation}\label{vareq}
\Hau^1(f(E))\leq K\,m^*(E),
\end{equation}
where $m^*$ is the outer Lebesgue measure.
\end{lemma}

\begin{proof} If $E$ is finite or denumerable, then the contidion~\eqref{vareq}
follows trivially. Suppose that $E$ is not denumerable. Let $\e>0$ be given, and
let $A$ be an open subset
of $[a,b]$ such that $E\subset A$ and $m(A)\leq m^*(E)+\e$.
Define inductively $E_0:=\emptyset$, and
\begin{equation*}
\begin{split}
E_i:=\{&x\in A\setminus E_{i-1}: B(x,1/i)\subset A\text{ and }\\
&\rho(f(x+t),f(x))\leq (K+\e)|t|\text{ for }|t|<1/i\}\ \text{ for }i\in\N.
\end{split}
\end{equation*}
Then each $E_i$ is Borel (see~\eqref{embed} in conjunction with e.g.\ \cite[Lemma~2.3]{D}).
Let $E_{ij}$ be such that $\diam(E_{ij})<1/i$, $(E_{ij})_j$ is
a pairwise-disjoint collection of Borel sets for each $i$, and $\bigcup_j E_{ij}=E_i$.
Note that $E\subset \bigcup_i E_i$.
We see that $f|_{E_{ij}}$ is $(K+\e)$-Lipschitz. 
It easily follows (see~\cite[Theorem~2.10.11]{F}) that
\[\Hau^1(f(E_{ij}))\leq (K+\e)\,m(E_{ij}),\]
and thus
\begin{equation*}
\begin{split}
\Hau^1(f(E))&\leq\Hau^1\bigg(f\bigg(\bigcup_{i,j}E_{ij}\bigg)\bigg)
\leq (K+\e)\sum_{i,j} m(E_{ij})\\
&\leq (K+\e)\,m(A)
\leq (K+\e)\,(m^*(E)+\e).
\end{split}
\end{equation*}
To obtain~\eqref{vareq}, send $\e\to0$.
\end{proof}

We have the following metric analogue of~\cite[Theorem~1]{Var}:

\begin{theorem}\label{vathm1}
Let $f:[a,b]\to (M,\rho)$ be arbitrary, and let $E$ be any
measurable set on which $md(f,\cdot)$ is finite. Then
\begin{equation}\label{vaeq2}
\Hau^1(f(E))\leq \int_E md(f,x)\,dx.
\end{equation}
\end{theorem}

\begin{proof} Using Lemma~\ref{funlem}, the proof is similar to the proof of~\cite[Theorem~1]{Var}.
Here are the details. First suppose that $md(f,x)<B$ for some $B\in\N$ on $E$. Let
\[ E_{nk}=\{x\in E:k-1\leq 2^n\cdot md(f,x)<k\},\ k=1,\dots, B2^n, n=1,\dots.\]
Then for each $n\in\N$ we have
\begin{equation*}
\begin{split}
\Hau^1 (f(E))&=\Hau^1 \bigg(f\bigg(\bigcup_k E_{nk}\bigg)\bigg)=\Hau^1\bigg(\bigcup_k f(E_{nk})\bigg)\leq \sum_k\Hau^1(f(E_{nk}))\\
&\leq\sum_k\frac{k}{2^n} m(E_{nk})
=\sum_k \frac{k-1}{2^n} m(E_{nk})+\frac{1}{2^n}\sum_k m(E_{nk}),
\end{split}
\end{equation*}
where the second inequality follows from Lemma~\ref{funlem}.
Therefore,
\[ \Hau^1(f(E))\leq \lim_{n\to\infty}\bigg[\sum_k \frac{k-1}{2^n} m(E_{nk})+\frac{1}{2^n}\sum_k m(E_{nk})\bigg]=\int_E md(f,x)\,dx.\]
\par
Now, if $md(f,x)$ is not bounded on $E$, then let 
\[ A_k=\{x\in E:k-1\leq md(f,x)<k\},\ \ k=1,\dots,\]
\begin{equation*}
\begin{split}
\Hau^1 (f(E))&=\Hau^1 \bigg(f\bigg(\bigcup_k A_k\bigg)\bigg)=\Hau^1 \bigg(\bigcup_k f(A_k)\bigg)\leq \sum_k\Hau^1(f(A_k))\\
&\leq \sum_k \int_{A_k} md(f,x)\,dx=\int_E md(f,x)\,dx.
\end{split}
\end{equation*}
\end{proof}

%

\section{Absolutely continuous functions}

We say that $f:[a,b]\to(M,\rho)$ has
{\em $($Luzin's$)$ property $(N)$} provided 
\begin{equation}\label{propertyN}
{\mathcal H}^1 (f(B))=0 \ \  \text{whenever}
\ \ B\subset[0,1]\text{ with }m(B)=0.
\end{equation}

The proof of the following theorem is standard (see e.g.\ \cite{S}
and the proof of Theorem in~\cite{DZ}).

\begin{theorem}\label{ACNthm} Let $f:[a,b]\to (M,\rho)$ is absolutely continuous.
Then $f$ has the property~$(N)$.
\end{theorem}

The previous theorem has the following corollary.

\begin{corollary} An absolutely continuous function $f:[a,b]\to(M,\rho)$ 
maps measurable subsets of $[a,b]$ onto $\Hau^1$-measurable subsets of $M$.
\end{corollary}

We will need the following theorem (see~\cite[Theorem~14]{Var} for the real-valued case).

\begin{theorem}\label{DZprop} 
Let $f:[a,b]\to(M,\rho)$ be continuous and has bounded variation.
Then $md(f,\cdot)$ exists almost everywhere in~$[a,b]$, is integrable, and
\begin{equation}\label{dzeq} 
\int^b_a md(f,x)\,dx\leq \bigvee^b_a f.
\end{equation}
\par
Further, if $f$  has the property~$(N)$, then the equality holds.
\end{theorem}

\begin{proof} Denote $A=[a,b]$. Theorem~\ref{Dthm} (ii) implies that $md(f,\cdot)$ 
exists for all $x\in A\setminus N$ with $m(N)=0$.
The area formula~\cite[Theorem~2.12]{D} together
with~\cite[Theorem~2.10.13]{F} implies that 
\begin{equation}\label{vypocet}
\bigvee^b_a f=\int N(f|_A,y)\,d\Hau^1 y\geq\int_{f(A\setminus N)}N(f|_A,y)\,d\Hau^1 y
=\int_{A\setminus N} md(f,x)\,dx,
\end{equation}
and thus~\eqref{dzeq} holds. Here, $N(f|_A,y)$ is the number of $x\in A$ such that $f(x)=y$.
\par
If $f$ has property~$(N)$, then clearly $\Hau^1(f(N))=0$ and we get equality
instead of an inequality in~\eqref{dzeq}. To see that, 
we have (again using~\cite[Theorem~2.12]{D} together
with~\cite[Theorem~2.10.13]{F})
\begin{equation*}
\begin{split}
\bigvee^b_a f&=\int N(f|_A,y)\,d\Hau^1 y=\int_{f(A)\setminus f(N)}N(f|_A,y)\,d\Hau^1 y\\
&\leq\int N(f|_{A\setminus N},y)\,d\Hau^1 y=\int_{A\setminus N} md(f,x)\,dx\leq\int_A md(f,x)\,dx.
\end{split}
\end{equation*}
\end{proof}

\begin{remark}\label{eqrem}
If $f$ from the previous theorem is absolutely continuous, then
we have equality in~\eqref{dzeq} (as $f$ has bounded variation, and it also
satisfies property~$(N)$ by Theorem~\ref{ACNthm}). It is easy to see
that $f$ is absolutely continuous if and only if $v_f$ is.
If that is the case, then
it follows that
\[ \int^d_c md(f,x)\,dx=v_f(d)-v_f(c)=\int^d_c v_f'(x)\,dx,\]
for each interval $[c,d]\subset[a,b]$. It is easy to see that
$md(f,x)\leq v_f'(x)$ whenever $v_f'(x)$ and $md(f,x)$ exist. 
Thus if $f$ is absolutely continuous, 
then $md(f,x)=v_f'(x)$ almost everywhere.
\end{remark}

The following version of the Banach-Zarecki theorem (see e.g.~\cite{N} or~\cite{Var} for the real-valued statement) 
was proved by L.\ Zaj\'\i\v{c}ek and the author in~\cite{DZ} using a result of Luzin~\cite{L} and a theorem about 
the Banach indicatrix function from~\cite{F}.
Here, we present a different proof, which is in the spirit of Varberg's approach (see~\cite[Theorem~3]{Var}).

\begin{theorem}\label{banzar}
Let $f:[a,b]\to(M,\rho)$ be a function. Then $f$ is absolutely continuous if and 
only if $f$ is continuous, has bounded variation, and has the property~$(N)$.
\end{theorem}

\begin{proof}
If $f$ is absolutely continuous, then a standard argument shows that
$f$ is continuous, and $f$ has bounded variation. Theorem~\ref{ACNthm}
shows that $f$ also has the property~$(N)$.
\par
To prove the converse, let $[a_i,b_i],i=1,\dots,k,$ be non-overlapping intervals in~$[a,b]$,
and let $E_i=\{x\in[a_i,b_i]:md(f,x)\text{ exists}\}$. Since by Theorem~\ref{Dthm} (ii)
we have that $m([a_i,b_i]\setminus E_i)=0$, and since $f$ has the property~$(N)$,
we obtain $\Hau^1\big(f(E_i)\big)=\Hau^1\big(f([a_i,b_i])\big)$.
Therefore,
\begin{equation}\label{altbz}
\begin{split}
\sum^k_{i=1} \rho(f(b_i),f(a_i))&\leq \sum^k_{i=1} \Hau^1(f([a_i,b_i]))=
\sum^k_{i=1} \Hau^1(f(E_i))\\
&\leq \sum^k_{i=1} \int_{E_i} md(f,x)\,dx 
=\sum^k_{i=1} \int^{b_i}_{a_i} md(f,x)\,dx,
\end{split}
\end{equation}
where the first inequality follows from~\cite[Corollary~2.10.12]{F} and the
second from Theorem~\ref{vathm1}. It is easy to see 
that the rightmost term in~\eqref{altbz} goes to~$0$, as $\sum^k_{i=1} (b_i-a_i)\to0$.
This last property follows from the fact that $md(f,\cdot)$ is integrable
by Theorem~\ref{DZprop}, and from a well-known property of the integral.
\end{proof}

The proofs of the next two theorems are analogous to the previous one
(cf.\ \cite[Theorems~4,~5]{Var}).

\begin{theorem}\label{Va4} If $f:[a,b]\to (M,\rho)$ is continuous, $md(f,\cdot)$
exists for all but finite or denumerable set of points and $md(f,\cdot)$
is integrable on~$[a,b]$, then $f$ is absolutely continuous on~$[a,b]$.
\end{theorem}

\begin{theorem}\label{Va5} If $f:[a,b]\to (M,\rho)$ is continuous, $md(f,\cdot)$
exists almost everywhere and is integrable on~$[a,b]$, 
and if $f$ has the property~$(N)$, then $f$ is absolutely continuous on~$[a,b]$.
\end{theorem}

The next theorem is a~consequence of Theorem~\ref{Va4} -- see~\cite[p.~266]{N}
or~\cite[Theorem~6]{Var} for the real-valued version.

\begin{theorem} If $md(f,x)$ exists for all $x\in[a,b]$, and if $md(f,x)$
is integrable, then $f$ is absolutely continuous on~$[a,b]$.
\end{theorem}

The following theorem is an analogue of~\cite[Theorem~30.12]{Vai}.

\begin{theorem} Let $f:[a,b]\to (M,\rho)$ be continuous.
Assume that 
\begin{enumerate}
\item there exists a closed and denumerable $E\subset[a,b]$ such that that $f$
is absolutely continuous on each closed interval in $[a,b]\setminus E$, and
\item $\int^b_a md(f,x)\,dx<\infty$. 
\end{enumerate}
Then $f$ is absolutely continuous on~$[a,b]$.
\end{theorem}

\begin{proof} We will prove that $f$ satisfies the assumptions of Theorem~\ref{Va5}.
By Theorem~\ref{Dthm} (ii) we have that $md(f,x)$ exists almost everywhere in~$[a,b]$,
and the integrability of~$md(f,\cdot)$ follows from~(ii).
Let $(a_i,b_i),\ (i\in\mcI\subset\N)$ be the intervals contiguous to $E$ in $[a,b]$.
%
By Theorem~\ref{ACNthm} and condition~(i), it follows that
$f|_{[a_i,b_i]}$ has property~$(N)$ for each $i\in\mcI$. As $E$ is denumerable,
we easily obtain that $f$ has property $(N)$. Thus Theorem~\ref{Va5} applies
and $f$ is absolutely continuous on~$[a,b]$.
\end{proof}

We have the following (see also~\cite[p.~246]{N} or~\cite[Theorem~9]{Var}):

\begin{theorem}
Let $f:[a,b]\to(M,\rho)$ be an absolutely continuous function, and $md(f,x)=0$ almost
everywhere on~$[a,b]$. Then $f$ is a constant function.
\end{theorem}

\begin{proof}
Theorem~\ref{DZprop} implies that $\bigvee^b_a f=0$. The only functions
with zero variation are the constant ones.
\end{proof}

The following theorem is an analogue of~\cite[Theorem~13]{Var}.

\begin{theorem}Let $f:[a,b]\to(M,\rho)$ be one-to-one and have bounded variation, let 
$A$ be any measurable set, and let
$E$ be the set of all $x\in A$ where $md(f,x)$ exists.
Then 
\begin{equation}\label{eqineq}
\int_A md(f,x)\,dx=\Hau^1(f(E))\leq \Hau^1(f(A)).
\end{equation}
The equality holds provided $f$ is absolutely continuous.
\end{theorem}

\begin{proof} 
First, assume that $f$ is absolutely continuous.
Then~\cite[Theorem~2.12]{D} shows that 
\begin{equation*} 
\begin{split}
\int_A md(f,x)\,dx&=\int_E md(f,x)\,dx
=\int N(f|_E,y)\,d\Hau^1y\\
&=\Hau^1(f(E))=\Hau^1(f(A)),
\end{split}
\end{equation*}
as $m(A\setminus E)=0$ by Theorem~\ref{Dthm} (ii), and $f$ has property~$(N)$
by Theorem~\ref{ACNthm}.
\par
Now, we will prove the equality from~\eqref{eqineq} for $f$, which are one-to-one
with bounded variation (note that the inequality in~\eqref{eqineq} holds trivially).
Define
\[A'_n:=\{x\in E:\rho(f(x+t),f(x))\leq n |t|\ \text{for }|t|<1/n\},\]
and $A_n:=A'_n\setminus_{j<n} A'_j$. Then each $A_n$ is
measurable (see e.g.\ \cite[Lemma~2.3]{D} together with~\eqref{embed}) and $A=\bigcup_n A_n$.
Further, write $A_n=\bigcup_k A_{nk}$ so that $(A_{nk})_k$ is
a pairwise-disjoint sequence of measurable sets with $\diam(A_{nk})<1/n$
for each $k$. Now extend each $f|_{A_{nk}}$ (which is $n$-Lipschitz by 
the definition of~$A_{nk}$) to a one-to-one $n$-Lipschitz function on $[a,b]$
(first extend $f|_{A_{nk}}$ to $\overline{A_{nk}}$ by continuity,
and then linearly and continuously on the intervals contiguous
to $\overline{A_{nk}}$; it is easy to see that the resulting function
is $n$-Lipschitz and one-to-one) --
call the extensions~$f_{nk}$.
Then
\begin{equation*}
\begin{split} 
\int_A md(f,x)\,dx&=\sum_{n,k}\int_{A_{nk}} md(f_{nk},x)\,dx=\sum_{n,k} \Hau^1(f_{nk}(A_{nk}))\\
&=\Hau^1\bigg(\bigcup_{n,k}f(A_{nk})\bigg)
=\Hau^1(f(E)),
\end{split}
\end{equation*}
where the first equality follows from the fact that almost all points of $A_{nk}$
are points of density, and $md(f,x)=md(f_{nk},x)$ at all such points
(see e.g.~\cite[Lemma~2.1]{D}).
The second equality follows by the previous paragraph.
\end{proof}

For the real-valued version of the following theorem, see~\cite[Theorem~15]{Var}.

\begin{theorem}\label{patnact} 
If $f:[a,b]\to(M,\rho)$ has bounded variation, and $A$
is a measurable subset of $[a,b]$, then $m^* (v_f(A))\geq \int_A md(f,x)\,dx$.
The equality holds if $f$ is absolutely continuous.
\end{theorem}

\begin{proof}
Let $E$ be the subset of $A$ where $md(f,x)$ exists.
Thus
\begin{equation}\label{va34}
\begin{split}
m^* (v_f(A))&\geq m^*(v_f(E))=\int_E v_f'(x)\,dx \\
&\geq\int_E md(f,x)\,dx=\int_A md(f,x)\,dx,
\end{split}
\end{equation}
where the first equality follows from~\cite[Theorem~13]{Var},
and the second inequality from the fact that $|v_f(y)-v_f(x)|\geq \rho(f(y),f(x))$
for all $x,y\in[a,b]$.
Note that we can write equalities instead of
inequalities in~\eqref{va34} provided $f$ is absolutely continuous
(as in that case Remark~\ref{eqrem} implies that $v_f$ is absolutely continuous,
and that $md(f,x)=v'_f(x)$ for almost every $x\in[a,b]$).
\end{proof}

Note that applying Lemma~\ref{funlem} with $K=0$ yields the following version of Sard's
theorem (cf.\ \cite[Lemma~2.2]{DZcurve}).

\begin{theorem}\label{Sard}
Let $f:[a,b]\to (M,\rho)$ and $E=\{x\in[a,b]:md(f,x)=0\}$. Then $\Hau^1(f(E))=0$.
\end{theorem}

The following is an analogue of~\cite[Theorem~18]{Var}.

\begin{theorem}\label{Va18}
Let $f:[a,b]\to(M,\rho)$ be continuous and of bounded variation.
Let $N$ be any set such that $\Hau^1(f(N))=0$. Then $m(v_f(N))=0$.
\end{theorem}

\begin{proof}
Denote $A=[a,b]$ and let $\e>0$. Let $K$ be a compact subset of $f(A)$
such that $K\cap f(N)=\emptyset$ and 
\[ \int_K N(f|_A,y)\,d\Hau^1y\geq \int N(f|_A,y)\,d\Hau^1y-\e;\]
existence of such a set $K$ follows from the regularity of $\Hau^1$
(see e.g.~\cite[\S 2.10.48]{F}).
Then $H=f^{-1}(K)$ satisfies $H\cap N=\emptyset$, and $v_f(H)\cap v_f(N)$
is at most denumerable. Now by~\cite[Theorem~2.10.13]{F} we have
\[ \bigvee^b_a f=\int N(f|_A,y)\,d\Hau^1y,\]
and thus
\begin{equation*}
\begin{split} 
m^*(v_f(H)) &\geq \int_H md(f,x)\,dx=\int N(f|_H,y)\,d\Hau^1y\\
&=\int_K N(f|_A,y)\,d\Hau^1y \geq \int N(f|_A,y)\,d\Hau^1y-\e\\
&= \bigvee^b_a f-\e,
\end{split}
\end{equation*}
where the first inequality follows from Theorem~\ref{patnact}, and the first
equality from~\cite[Theorem~2.12]{D}.
It follows that $m^*(v_f(N))\leq \e$, and as $\e>0$ was arbitrary, we have
that $m(v_f(N))=0$.
\end{proof}

The following is an analogue of~\cite[Theorem~19]{Var}.

\begin{theorem} Let $f:[a,b]\to (M,\rho)$ be continuous, have bounded variation, and let 
$E$ be a measurable set for which $\Hau^1(f(E))=0$. Then $md(f,x)=0$
for almost all $x\in E$.
\end{theorem}

\begin{proof} We have
\[ \int_E md(f,x) \leq m(v_f(E))\leq 0,\]
where the first inequality follows from Theorem~\ref{patnact}, and the second from
Theorem~\ref{Va18}.
\end{proof}

The following theorem was established by Vall\'ee Poussin~\cite{PO} for real
valued functions.

\begin{theorem}\label{pothm}
Let $(M,\rho)$ be a metric space and $f:[c,d]\to M$, $g:[a,b]\to[c,d]$ be absolutely
continuous functions. Then $f\circ g$ is absolutely continuous if and only if
$md(f,g(x))\cdot g'(x)$ is integrable.
\end{theorem}

\begin{remark}
The expression $h(x)=md(f,g(x))\cdot g'(x)$ is interpreted in the following
sense (usual in the measure-theory): $h(x)=0$ provided $g'(x)=0$ (even
when $md(f,g(x))$ does not exist).
\end{remark}

\begin{proof} Suppose that $f\circ g$ is absolutely continuous. Then Theorem~\ref{DZprop}
implies that $md(f\circ g,x)$ is integrable.
Let $A$ be the set of all points $x$ of $[a,b]$ where
$g(x)\neq0$ and $md(f\circ g,x)$ exists. 
Lemma~\ref{diflem} shows that for every $x\in A$, the metric derivative
$md(f,g(x))$ exists.
Thus, if $x\in A$, then we have that
$md(f\circ g,x)=md(f,g(x))\cdot g'(x)$ by a chain rule
for metric derivatives (see e.g.\ \cite[Lemma~2.4 (ii)]{DZcurve}).
Let $N:=\{x\in [a,b]:g'(x)=0\}$.
Then $m([a,b]\setminus(A\cup N))=0$ (by Theorem~\ref{Dthm} (ii)), and
thus $md(f,g(x))\cdot g'(x)$ is integrable on $[a,b]$.
\par
Suppose that $md(f,g(x))\cdot g'(x)$ is integrable. It is easily seen
that $f\circ g$ has property $(N)$ (as it is stable under compositions), 
and thus by Theorem~\ref{banzar} it is enough to show that $f\circ g$
has bounded variation. 
Let 
\[ A:=\{x\in[a,b]:md(f,g(x))\text{ exists and }g'(x)\neq0\},\]
$B_1:=\{x\in[a,b]:g'(x)=0\}$, and $B_2:=[a,b]\setminus (A\cup B_1)$.
Note that for almost every $x\in[a,b]$ we have that either
$g'(x)=0$ or $md(f,g(x))$ and $g'(x)\neq0$ exist (in the second case,
we also have that $md(f\circ g,x)$ exists and is equal to $md(f,g(x))\cdot g'(x)$
by the chain rule for metric derivatives~\cite[Lemma~2.4 (ii)]{DZcurve}). Thus
it follows that $m(B_2)=0$.
Let $B:=B_1\cup B_2$. By Theorem~\ref{Sard}, and because $f\circ g$
has property $(N)$, we have that $\Hau^1((f\circ g)(B))=0$.
We obtain
\begin{equation*}
\begin{split} 
\bigvee^b_a (f\circ g)&=\int_M N(f\circ g,y)\,d\Hau^1y
=\int_{M\setminus (f\circ g)(B)} N(f\circ g,y)\,d\Hau^1y\\
&=\int_{\{x\in[a,b]:f\circ g(x)\not\in B\}} md(f\circ g,x)\,dx\\
&\leq\int_A md(f,g(x))\cdot g'(x)\,dx<\infty,
\end{split}
\end{equation*}
where the first equality 
follows from~\cite[Theorem~2.10.13]{F}, and the third by~\cite[Theorem~2.12]{D}.
We have that $f\circ g$ has finite variation, and thus we can apply Theorem~\ref{banzar}.
\end{proof}


\begin{thebibliography}{WWW}



\bibitem{BL}
Y.~Benyamini, J.~Lindenstrauss, {\em Geometric Nonlinear Functional
Analysis, Vol.~1,} Colloquium Publications \textbf{48},
American Mathematical Society, Providence, 2000.

\bibitem{D}
J.~Duda, {\em Metric and $w^*$-differentiability of pointwise Lipschitz mappings},
submitted
(available electronically at \texttt{http://www.karlin.mff.cuni.cz/kma-preprints}).

\bibitem{DZcurve}
J.~Duda, L.~Zaj\'\i\v{c}ek, {\em Curves with values in Banach spaces --
differentiability via homeomorphisms}, to appear in the Rocky Mountain
J. of Math.


\bibitem{DZ}
J.~Duda, L.~Zaj\'\i\v{c}ek, {\em The Banach-Zarecki theorem for functions
with values in metric spaces}, to appear in the Proc. Amer. Math. Soc.

\bibitem{F}
H.~Federer, {\em Geometric Measure Theory,} Grundlehren der math. Wiss., vol.~153, Springer, New York, 1969.

\bibitem{K}
B.~Kirchheim, {\em Rectifiable metric spaces: local structure and
regularity of the Hausdorff measure}, Proc. Amer. Math. Soc.
~\textbf{121} (1994), 113--123.

\bibitem{L}
N.~Luzin, {\em The Integral and Trigonometric Series} (Russian),
 Mat. Sbornik \textbf{30} (1916), 1--242. 

\bibitem{N}
I.~P.~Natanson, {\em Theory of functions of a real variable} (Translated from Russian),
Revised Edition, Ungar, New York, 1961.

\bibitem{PO}
Ch.~J. de la Vall\'ee Poussin, {\em Sur l'integrale de Lebesgue}, Trans. Amer. Math. Soc.,
\textbf{16} (1915), 435--501.

\bibitem{S}
S.~Saks, {\em Theory of the integral}, Monographie Mat., vol. 7, Hafner, New York, 1937. 

\bibitem{Var}
D.~E.~Varberg, {\em
On absolutely continuous functions},
Amer. Math. Monthly \textbf{72} (1965), 831--841.

\bibitem{Vai}
J.~V\"ais\"al\"a, {\em Lectures on $n$-dimensional quasiconformal mappings}, 
Lecture Notes in Mathematics, vol. 229, Springer, Berlin-New York, 1971. 
 






\end{thebibliography}
\end{document}